\documentclass{statsoc}
\usepackage{amsthm,amsmath,amssymb} % amsfonts
\usepackage[a4paper]{geometry}
\usepackage[T1]{fontenc}
\usepackage[utf8]{inputenc}
\usepackage{natbib}
\usepackage{mathrsfs}  % mathscr
\usepackage[colorlinks,breaklinks,linkcolor=black,urlcolor=black,citecolor=black]{hyperref} % ,pagebackref

\title[Nonparametric independence tests in metric spaces]{Nonparametric independence tests in metric spaces: What is known and what is not}
\author[Castro-Prado and Gonz\'alez-Manteiga]{Fernando Castro-Prado}
\address{University of Santiago de Compostela and Health Research Institute, Santiago de Compostela, Spain.}
\email{fernando.castro.prado@rai.usc.es}
\author[Castro-Prado and Gonz\'alez-Manteiga]{Wenceslao Gonz\'alez-Manteiga}
\address{University of Santiago de Compostela, Santiago de Compostela, Spain.}

\DeclareMathOperator{\dcov}{dcov}
\DeclareMathOperator{\dvar}{dvar}
\DeclareMathOperator{\dcor}{dcor}
\DeclareMathOperator{\dCov}{dCov}

\DeclareMathOperator{\dCor}{dCor}
\newcommand{\dcovh}{\widehat{\dcov}}

\DeclareMathOperator{\E}{E}
\DeclareMathOperator{\Prob}{P}

\DeclareMathOperator{\Normal}{N}

\newcommand{\implica}{\Rightarrow}

\newcommand{\eqv}{\Leftrightarrow}
\newcommand{\flecha}{\longrightarrow}
\newcommand{\flechita}{\mapsto}

\newcommand{\inner}[2]{\left\langle{#1},{#2}\right\rangle}
\newcommand{\norma}[1]{\left\lVert#1\right\rVert}
\newcommand{\llaves}[1]{\left\lbrace{#1}\right\rbrace}
\newcommand{\recuadro}{\framebox[1.1\width]}

\newcommand{\overeq}[1]{\overset{\text{#1}}{=}}
\newcommand{\overdef}[1]{\overset{\text{#1}}{:=}}
\newcommand{\overleq}[1]{\overset{\text{#1}}{\leq}}
\newcommand{\overless}[1]{\overset{\text{#1}}{<}}

\newcommand{\overthen}[1]{\overset{\text{#1}}{\implica}}

\newcommand{\overarrow}[1]{\overset{\text{#1}}{\flecha}}
\newcommand{\detlim}{\underset{n\to\infty}{\flecha}}

\newcommand{\aslim}{\overset{a.s.}{\detlim}}
\newcommand{\distrilim}{\overset{\mathcal{D}}{\detlim}}

\newcommand{\disjoint}{\sqcup}
\newcommand{\comp}{\circ}
\newcommand{\card}{\#}
\newcommand{\tensor}{\otimes}

\newcommand{\wrt}{\,\mathrm{d}}
\newcommand{\dmu}{\,\mathrm{d}\mu}
\newcommand{\dnu}{\,\mathrm{d}\nu}
\newcommand{\damu}{\,\mathrm{d}|\mu|}
\newcommand{\amu}{a_\mu}
\newcommand{\anu}{a_\nu}
\newcommand{\dtheta}{\,\mathrm{d}\theta}
\newcommand{\LyonX}{\,\mathscr{X}}
\newcommand{\LyonY}{\,\mathscr{Y}}
\newcommand{\spx}{\mathscr{X}}
\newcommand{\spy}{\mathscr{Y}}
\newcommand{\spz}{\mathscr{Z}}
\newcommand{\xy}{\mathscr{X}\times\mathscr{Y}}
\newcommand{\dx}{d_{\mathscr{X}}}
\newcommand{\dy}{d_{\mathscr{Y}}}
\newcommand{\dz}{d_{\mathscr{Z}}}
\newcommand{\xd}{\left(\mathscr{X},d_{\mathscr{X}}\right)}
\newcommand{\yd}{\left(\mathscr{Y},d_{\mathscr{Y}}\right)}
\newcommand{\xsd}{\left(\mathscr{X},\sqrt{d_{\mathscr{X}}}\right)}
\newcommand{\ysd}{\left(\mathscr{Y},\sqrt{d_{\mathscr{Y}}}\right)}
\newcommand{\intx}{\int_{\mathscr{X}}}
\newcommand{\h}{\mathcal{H}}
\newcommand{\Borel}[1]{\mathcal{B}\left({#1}\right)}
\newcommand{\iid}{\text{ i.i.d. }}

\newcommand{\R}{\mathbb{R}}
\newcommand{\Rplus}{\mathbb{R}^{+}}
\newcommand{\Q}{\mathbb{Q}}
\newcommand{\Z}{\mathbb{Z}}
\newcommand{\N}{\mathbb{N}}
\newcommand{\Nstar}{\mathbb{N}^{*}}
\newcommand{\Zplus}{\mathbb{Z}^{+}}

\begin{document}

\begin{abstract}
Distance correlation is a recent extension of Pearson's correlation, that characterises general statistical independence between Euclidean-space-valued random variables, not only linear relations. This review delves into how and when distance correlation can be extended to metric spaces, combining the information that is available in the literature with some original remarks and proofs, in a way that is comprehensible for any mathematical statistician.
\end{abstract}

\keywords{Distance correlation; Association measures; Nonparametric statistics}

\section{Introduction}
The \emph{energy of data} \citep{TEOD} and all the mathematical statistics that stems from it, including the characterisation of independence in Euclidean spaces (\S~\ref{SRB}) and many other interesting results \citep{DISCO,Brownian,JMultivar}, have a very strong and well-established theoretical basis \citep{Bakirov,SRB,TEOD}.

Nevertheless, the article \citep{Lyons} that introduces distance correlation in metric spaces leaves a surprising amount of details to the reader \citep[p. 2]{Jakobsen}. The elision of so many intermediate steps meant that, for several years, it was unnoticed that most of the theory was incorrect \citep{Ups}. Such mistakes were largely discovered by \citet{Jakobsen}, who devoted 150 pages to go through and to correct glitches of the original 10-page paper.

The goal of the present review is to present a corrected version of Lyons' theory, by summarising and explaining the work by \citet{Jakobsen} and by adding a few original proofs, all of this taking into account the recent corrigendum of the original article \citep{Ups}. In addition, the reader will be provided with a gentle introduction to the abstract mathematical concepts that this theory requires. Thus, for the first time, a clear and concise bottom-up explanation of the theory of distance correlation in metric spaces is available to the scientific community.

\section{Distance correlation in Euclidean spaces}\label{SRB}
When two random elements (vectors) $X$ and $Y$ are Euclidean-space-valued (let $X$ be $L-$dimensional and $Y$ be $M-$dimensional; for $L,M\in\Zplus$), it is possible to define an association measure that characterises their independence called \emph{distance correlation} \citep{SRB}. Firstly, distance covariance should be defined, as a certain norm of the difference of the joint characteristic function and the product of the marginals:
\[\dCov(X,Y):=\norma{\varphi_{X,Y}-\varphi_X\varphi_Y}_w\equiv\sqrt{\int_{\mathbb{R}^L\times\mathbb{R}^M}|\varphi_{X,Y}(t,s)-\varphi_X(t)\varphi_Y(s)|^2w(t,s)\,\mathrm{d}t\,\mathrm{d}s}\text{;}\]\setcitestyle{square}%
where $w$ is a weight function which is dependent of the dimension of the Euclidean spaces in which the supports of $X$ and $Y$ are contained (and it has a property of uniqueness \citep{Uniqueness}):\setcitestyle{round}%
\[
w(t,s):=
\frac{\Gamma\left(\frac{L+1}{2}\right)}{\left(\norma{t}\sqrt{\pi}\right)^{L+1}}\:
\frac{\Gamma\left(\frac{M+1}{2}\right)}{\left(\norma{s}\sqrt{\pi}\right)^{M+1}}
,\;(t,s)\in\mathbb{R}^L\times\mathbb{R}^M\text{.}
\]
And, as usually:
\[\varphi_X(t):=\E\left[e^{i\inner{t}{X}}\right],\:t\in\mathbb{R}^L; \;\;\;\;\;\;
\varphi_Y(s):=\E\left[e^{i\inner{s}{Y}}\right],\:s\in\mathbb{R}^M\text{.}\]

Logically, distance correlation is defined as the quotient of variance and the product of standard deviations and so it has no sign:
\[\dCor(X,Y):=\frac{\dCov(X,Y)}{\sqrt{\dCov(X,X)\dCov(Y,Y)}}\text{,}\]
whenever $\dCov(X,X)\dCov(Y,Y)\neq0$. If $\dCov(X,X)=0$, then $\dCor(X,Y):=0$.

The reasons why distance correlation is an improved version of the squared (Pearson's) correlation are:
\begin{itemize}
	\item It has values in [0,1]. This is unsurprising: $\R$ is totally ordered and, as such, one can only move ``leftwards'' or ``rightwards'' and so the sign of (Pearson's) correlation expresses this structure. However, this notion is not valid in Euclidean spaces of arbitrary dimensionality.
	\item It is zero if \underline{and only if} $X$ and $Y$ are independent (thus, its interest).
\end{itemize}
Notwithstanding the convoluted initial definition of \emph{dCor}, its sample version can easily be computed. Given a paired sample
\[(X_1,Y_1),\ldots,(X_n,Y_n)\iid(X,Y)\text{;}\]
let $a_{ij}:=d(X_i,X_j)$ for $i,j\in[1,n]\cap\mathbb{Z}$. Using this notation, doubly-centred distances are:
\[A_{ij}:=a_{ij}-\bar a_{i\cdot}-\bar a_{j\cdot}+\bar a_{\cdot\cdot}\]
If $\{b_{ij}\}_{i,j}$ and $\{B_{ij}\}_{i,j}$ are analogously defined for $\{Y_i\}_i$, the empirical distance covariance is simply the nonnegative real number whose square is:
\[\widehat{\dCov}_n(X,Y)^2:=\frac{1}{n^2}\sum_{i,j=1}^n A_{ij}B_{ij}\text{,}\]
so that it is, indeed, a correlation of distances.

The above estimator comes from the alternative definition of \textit{dCov} derived by \citet{Brownian}:
\[\dcov(X,Y)\!=\!\E[d(X,X')d(Y,Y')]+\E[d(X,X')]\E[d(Y,Y')]-2\E[d(X,X')d(Y,Y'')]\text{,}\]
which is valid as long as moments of order $2$ are finite. Primed letters refer to independent and identically distributed copies of the corresponding random element.

Whenever $\{X,Y\}$ are independent and have finite first moments, the asymptotic distribution of the product of a scaled version of the preceding statistic is a linear combination of independent chi-squared variables with one degree of freedom. More precisely:
\[n\:\widehat{\dCov}_n(X,Y)^2\distrilim\sum_{j=1}^\infty \lambda_j Z_j^2\text{,}\] 
where $\{Z_j\}_j$ are i.i.d. $\Normal(0,1)$ and $\{\lambda_j\}_j\subset\mathbb{R}$. Unfortunately, this null distribution is not useful in practice.

Instead, it is resampling techniques that should be used. The most sensible choice when it comes to approximating the null distribution of the test statistic is to base the design of the resampling scheme on the information that $H_0$ provides, which in this case (i.e., independence) yields to \emph{permutation tests}.

\section{Context and notations}\label{Intro_Jakobsen}
\subsection{General statement of the nonparametric problem of independence}
Let $\xd$ and $\yd$ be two arbitrary separable metric spaces (the need for separability is dealt with in~\ref{Sep}). The random element $Z=(X,Y)$ is defined over $(\Omega,\mathcal F,\Prob)$ and has values in $\spx\times\spy$, with its distribution being
\[\theta:\Borel{\spx\times\spy}\longrightarrow[0,1]\text{.}\]
The following notation will be used for the marginal distributions:
\begin{itemize}
	\item $X\sim\mu:=\theta\comp\pi_1^{-1}$, marginal over $\spx$; where $\pi_1:(x,y)\in\spx\times\spy\flechita x\in\spx$.
	\item $Y\sim\nu:=\theta\comp\pi_2^{-1}$, \hspace*{.5mm}marginal over $\spy$; \hspace*{.5mm}where $\pi_2:(x,y)\in\spx\times\spy\flechita y\in\spy$.
\end{itemize}
\begin{sloppypar}
Thus, the nonparametric test of independence for $X$ and $Y$ consists in testing ${H_0:\theta=\mu\times\nu}$ versus ${H_1:\theta\neq\mu\times\nu}$. For the sake of clarity, it is important to note that the product $\mu\times\nu$ is defined conventionally: it is the only measure in $\Borel{\spx}\tensor\Borel{\spy}$ so that
\end{sloppypar}
\[(\mu\times\nu)(A\times B):=\mu(A)\nu(B);\;A\in\Borel{\spx},\:B\in\Borel{\spy}\text{.}\]

\subsection{Separability of marginal spaces}\label{Sep}
The first perquisite of assuming the separability of $\spx$ and $\spy$ is that, this way, the $\sigma-$algebra that their topological product generates is simply the product $\sigma-$algebra:
\[\Borel{\spx\times\spy}=\Borel{\spx}\tensor\Borel{\spy}:=\sigma\llaves{A\times B:A\in\Borel{\spx},B\in\Borel{\spy}}\text{.}\]\setcitestyle{square}%

\begin{sloppypar}
This equality is useful by itself (e.g., it is crucial to the proof of lemma 3.10 in \citet{Jakobsen}), but its most important corollary is that it guarantees that the metrics of the marginal spaces are jointly measurable: for $\spz\in\llaves{\spx,\spy}$, $\dz$ is {$\Borel{\spz}\tensor\Borel{\spz}/\Borel{\mathbb{R}}-$measurable}. This, in turn, is what ensures that the Lebesgue integrals that appear in the definition of distance covariance (\S~\ref{Def_dcov}) are defined. A counterexample would be $\spx:=\mathbb{R}^{\mathbb{R}}$, equipped with the discrete metric. This is a particular case of \emph{Nedoma's pathology}
(see \citet[proposition  21.8]{Schechter} and \citet[example 6.4.3]{Bogachev} for further details), which states that the diagonal set $\{(x,x):x\in\spx\}$ is not in $\Borel{\spx}\tensor\Borel{\spx}$ when the cardinality of $\spx$ is greater than that of the continuum.
\end{sloppypar}\setcitestyle{round}%
	
Finally, separability is explicitly used in the proofs of some important properties of distance covariance \citep[theorem 4.4 and lemma 5.8]{Jakobsen}, which indicates that it is not an ungodly hypothesis.

The original article that presented distance correlation in metric spaces \citep{Lyons} was oblivious of the crucial role of separability in the theory.

\subsection{Signed measures}\label{Signed}
The map $\mu:\Borel{\spx}\longrightarrow\mathbb{R}$ is said to be a finite signed (Borel) measure, and it is denoted $\mu\in M(\spx)$, if and only if $|\mu|$ is a finite measure. For each $\mu\in M(\spx)$, there is a \emph{Hahn--Jordan decomposition} and it is essentially unique \citep[theorem 3.2.1]{Billingsley} or, in other words, it is possible to find a couple of nonnegative measures $\mu^{\pm}\in M(\spx)$ so that
\[\mu=\mu^+-\mu^-\]
and a partition of the space $\spx=\spx^+\disjoint\spx^-$ satisfying:
\[\mu^+(\spx^-)=0=\mu^-(\spx^+)\text{;}\]
which is to say that $\mu^+$ and $\mu^-$ are orthogonal (mutually singular).

This allows to naturally define (Lebesgue) integrals with respect to signed measures. For $f:\spx\longrightarrow\mathbb{R}$ measurable,
\[\intx f\dmu:=\intx f\dmu^+-\intx f\dmu^-\text{;}\]
which is well-defined whenever $f$ is integrable with respect to $|\mu|=\mu^++\mu^-$.

On the other hand, it will also be necessary to integrate with respect to product measures. To begin with, consider $\nu\in M(\spy)$, with Hahn--Jordan decomposition given by $(\spy^{\pm},\nu^{\pm})$. Then:
\begin{itemize}
	\item $\mu^+\times\nu^++\mu^-\times\nu^-$ is a (nonnegative) measure with support $(\spx^+\times\spy^+)\disjoint(\spx^-\times\spy^-)$;
	\item $\mu^+\times\nu^-+\mu^-\times\nu^+$ is a (nonnegative) measure with support $(\spx^+\times\spy^-)\disjoint(\spx^-\times\spy^+)$.
\end{itemize}
Because of their disjoint supports, the aforementioned two measures are mutually singular and, consequently \citep[corollary of theorem 6.14]{Rudin}, they form the Hahn--Jordan decomposition of $\mu\times\nu$:
\[\mu\times\nu=(\mu^+\times\nu^++\mu^-\times\nu^-)-(\mu^+\times\nu^-+\mu^-\times\nu^+)\text{.}\]
Thus, the integral of a Borel-measurable function $h:\spx\times\spy\longrightarrow\mathbb{R}$ with respect to $\mu\times\nu$ is:
\[\int h\wrt\mu\times\nu=\int h\wrt\mu^+\times\nu^++\int h\wrt\mu^-\times\nu^--\int h\wrt\mu^+\times\nu^--\int h\wrt\mu^-\times\nu^+\text{;}\]
which entails that $\mathcal L^1(\mu\times\nu)$ is the intersection of the four function spaces $\mathcal L^1(\mu^{\pm}\times\nu^{\pm})$.

On the last equation, the integration sets were omitted, as it is superfluous to underscore that it is the largest possible one (in this case, $\spx\times\spy$). This notation abuse, taken from \citet{Lyons}, is among the few ones that will be used on the present paper, while the ones that caused mistakes and confusion on Lyons' article (and even in its corrigendum \setcitestyle{square}\citep{Ups}) will be avoided. \setcitestyle{round}

The last relevant remark about the integration with respect to the product of signed measures is that they satisfy a generalised Fubini--Tonelli theorem \citep[\S~3.3]{Bogachev}:
\[\forall\:h\in\mathcal L^1(\mu\times\nu),\;\int h\wrt\mu\times\nu=\iint h\dmu\dnu=\iint h\dnu\dmu\text{.}\]

\subsection{Regularity of a measure}\label{Moments}

For the sake of clarity, it is convenient to state and prove the $c_r-$\emph{inequality}. For any $\alpha,\beta,r\in\Rplus$: $(\alpha+\beta)^r\leq c_r(\alpha^r+\beta^r)$, where
\[c_r=\begin{cases}1,&r<1\\2^{r-1},&r\geq1\end{cases}\text{.}\]
\textbf{Proof}. (1) Let \recuadro{r<1}. The goal is to show that 
\[(t+1)^r\leq t^r+1,\;t:=\frac{\alpha}{\beta}\]
or, equivalently, that
\[f(t):=t^r+1-(t+1)^r\geq0\text{.}\]
And the latter inequality holds because $r-1<0$:
\[\forall t\in\Rplus,\;f'(t)=r(t^{r-1}-(t+1)^{r-1})>0\implica\forall t\in\Rplus,\;f(t)\geq f(0)=0\text{.}\]
(2) For \recuadro{$r\geq1$}, the function $g(x):=x^r$ is convex in every $x\in\Rplus$. When $r>1$:
\[g''(x)=r(r-1)x^{r-2}>0,\;x\in\Rplus\text{.}\]
Geometrically, convexity implies that:
\[\pushQED{\qed} 
g\left(\frac{\alpha+\beta}{2}\right)\leq\frac{g(\alpha)+g(\beta)}{2}\eqv(\alpha+\beta)^r\leq2^{r-1}(\alpha^r+\beta^r)\text{.}\qedhere
\popQED
\]

At this point, it is possible to introduce the concept of regularity of a signed measure: $\mu\in M(\spx)$ is said to have finite moments of order $r$, and it is written as $\mu\in M^r(\spx)$, if and only if
\[\exists\:o\in\spx,\;\int\dx(o,x)^r\damu(x)<+\infty\text{.}\]
Applying the $c_r-$inequality, it is straightforward to see that when the condition above holds, it does so for any origin:
\[\mu\in M^r(\spx)\eqv\forall o\in\spx,\;\int\dx(o,x)^r\damu(x)<+\infty\text{.}\]
In addition, a signed measure on a product of two spaces $\theta\in M(\spx\times\spy)$ is said to belong to $M^{r,r}(\spx\times\spy)$ if both its marginals have finite moments of order $r$. Finally, the subindex $1$ will be used as a notation for probability measures:
\[M_1(\spx):=\big\lbrace{\mu\in M(\spx):\:\mu\geq0,\:\mu(\spx)=1}\big\rbrace;\]
\[M_1^r(\spx):=M^r(\spx)\cap M_1(\spx);\;\;\;\;M_1^{r,r}(\spx\times\spy):=M^{r,r}(\spx\times\spy)\cap M_1(\xy)\text{.}\]

\section{Formal definition of \emph{dcov}}\label{Def_dcov}
The previous section set the theoretical framework in which speaking of distance covariance makes sense, thus solving some inconsistencies of \citet{Lyons}. This will enable to define the operator \emph{dcov} rigorously, simplifying and illustrating the explanations by \citet{Jakobsen}.

\subsection{Integrability of the metric}\label{Int}
In order to define \emph{dcov}, it is important to keep in mind that:
\[\forall\:\mu_1,\mu_2\in M^1(\spx):\;\dx\in\mathcal L^1(\mu_1\times\mu_2)\text{.}\]
This is a consequence of Fubini and the triangle inequality:
\[\int\dx\wrt|\mu_1|\times|\mu_2|\leq\int\dx(x,o)\wrt|\mu_1|\times|\mu_2|(x,x')+\int\dx(o,x')\wrt|\mu_1|\times|\mu_2|(x,x')=\]
\[=|\mu_2|(\spx)\int\dx(x,o)\wrt|\mu_1|(x)+|\mu_1|(\spx)\int\dx(x,o)\wrt|\mu_2|(x)<+\infty\text{.}\]

\subsection{Expected distances and some inequalities}\label{aD}
The definition of distance covariance involves doubly centred distances (\S~\ref{dmu}), but first the various expected values that are to appear should be checked to be well-defined. For $\mu\in M^1(\spx)$, the following function maps each point $x\in\LyonX$ to its expected distance to the random element $X\sim\mu$:
\begin{align*}
\amu:\;&\spx\longrightarrow\mathbb{R}\\
&x\longmapsto\int\dx(x,x')\dmu(x')
\end{align*}
Obviously, it is well-defined. On top of that, it is $|\mu|(\spx)-$Lipschitzian  (and, therefore, continuous):
\[\forall x,x'\in\spx:\;|\amu(x)-\amu(x')|\leq\int|\dx(x,z)-\dx(x',z)|\damu(z)\leq\]
\[\leq\int\dx(x,x')\damu(z)=|\mu|(\spx)\dx(x,x')\text{.}\]
On the other hand, recalling~\ref{Int}, the integral $D(\mu)$ is always a real number:
\[D(\mu):=\int\amu\dmu=\int\dx\dmu\times\mu\text{.}\]
The following four inequalities can easily be derived from the previous results and they will be very useful hereinafter. For $\mu\in M_1^1(\spx)$ and $x,y\in\spx$:
\begin{enumerate}
	\item $D(\mu)\leq2\amu(x)$;
	\item $D(\mu)\leq\amu(x)+\amu(y)$;
	\item $\dx(x,y)\leq\amu(x)+\amu(y)$;
	\item $\amu(x)\leq\dx(x,y)+\amu(y)$.
\end{enumerate}
\textbf{Proof}.
\begin{sloppypar}
(1) $D(\mu)=\int\dx(x',x'')\dmu^2(x',x'')\leq$

$\leq\mu(\spx)\int\dx(x',x)\dmu(x')+\mu(\spx)\int\dx(x,x'')\dmu(x'')=2\amu(x)$.

(2) Applying \emph{(1)} to $x$ and $y$ and adding side-by-side the resulting equations, one gets: $2D(\mu)\leq2\amu(x)+2\amu(y)$.

(3) Integrate with respect to $\mu(z)$ both sides of: $\dx(x,y)\leq\dx(x,z)+\dx(y,z)$.

(4) Idem to \emph{(3)}: $\dx(x,z)\leq\dx(x,y)+\dx(y,z)$.\qed
\end{sloppypar}

\subsection{Doubly centred distances}\label{dmu}
For $\mu\in M^1(\spx)$, the doubly $\mu-$centred version of $\dx$ is:
\begin{align*}
d_\mu:\;&\spx\times\spx\longrightarrow\mathbb{R}\\
&(x_1,x_2)\flechita\dx(x_1,x_2)-\amu(x_1)-\amu(x_2)+D(\mu)
\end{align*}

This modification of $\dx$, in general, is not a metric; although it is always continuous (since $\dx$, $\amu$, $\pi_1$ and $\pi_2$ are) and, in particular, Borel-measurable. Moreover, it is important to note that, when writing $d_\mu$, there is no explicit reference to the metric space over which this map is defined. Such an abuse of notation makes formulae easier to read and write without creating any misunderstanding. That is not the case of some abbreviations by Lyons, such as the usage of $d:=\dx$ and $d:=\dy$, which mistakenly suggests that there is a need for $\LyonX$ and $\LyonY$ to share the same metric structure, which is an unnecessary restriction for the theory that would render some interesting applications impossible.

The last remarkable property of $d_\mu$ is:
\[\forall\mu,\mu_1,\mu_2\in M^1_1(\spx):\;d_\mu\in\mathcal L^2(\mu_1\times\mu_2)\text{.}\]

\textbf{Proof}. In the first instance, it is convenient to justify that, for any $(x,y)\in\spx^2$,
\[|d_\mu(x,y)|\leq2\amu(y)\text{.}\]
To see this, there are two cases to be considered:
\begin{itemize}
	\item{If \recuadro{$d_\mu(x,y)\geq0$}, it suffices to apply the inequalities in~\ref{aD}:
		\[|d_\mu(x,y)|=d_\mu(x,y)\overleq{(3)}D(\mu)\overleq{(1)}2\amu(y)\text{.}\]}
	\item{For \recuadro{$d_\mu(x,y)<0$}, the arguments of \citet[páx. 10]{Jakobsen} make use of unnecessarily strong hypotheses. Instead, the following rationale:
		\[\forall z,t\in\spx:\:\dx(x,z)\leq\dx(x,y)+\dx(y,t)+\dx(t,z)\implica
		\]
		\[\implica
		\amu(x)\leq\dx(x,y)+\amu(y)+D(\mu)\text{;}\]
		yields $|d_\mu(x,y)|\leq2\amu(y)$.}
\end{itemize}
Now, using the aforementioned inequality, proving that $d_\mu\in\mathcal L^2(\mu_1\times\mu_2)$ turns out to be quite straightforward:
\[\int d_\mu(x,y)^2\dmu_1\times\mu_2(x,y)\leq4\int\amu(x)\amu(y)\dmu_1\times\mu_2(x,y)\overeq{Fubini}\]
\[\pushQED{\qed}
=4\int\dx(x,z)\dmu_1\times\mu(x,z)\:\int\dx(y,z)\dmu_2\times\mu(y,z)\overless{\recuadro{$\dx\in\mathcal L^1$}}+\infty\text{.}\qedhere
\popQED\]

\subsection{The association measure \emph{dcov}}\label{dcov}
The generalised distance covariance is defined as:
\[\dcov(\theta):=\int_{(\spx\times\spy)^2}d_\mu(x,x')d_\nu(y,y')\dtheta^2\left((x,y),(x',y')\right),\;\theta\in M_1^{1,1}(\spx\times\spy)\text{;}\]
where, once again, $\mu:=\theta\comp\pi_1^{-1}$ and $\nu:=\theta\comp\pi_2^{-1}$.

In order to check that \emph{dcov} is well-defined, it suffices to note that the integral of the product of two functions with respect to a (nonnegative) measure is always a scalar product (bilinear, semidefinite positive) and, as a result, it satisfies the Cauchy--Bunyakovsky--Schwarz inequality. It is also possible to prove this particular case of H\"older's inequality more directly:
\[0\leq\int[d_\mu(v)d_\nu(w)-d_\mu(w)d_\nu(v)]^2\dtheta^2(v,w)=2\int d_\mu^2\dtheta^2\:\int d_\nu^2\dtheta^2-2\left(\int d_\mu d_\nu\dtheta^2\right)^2\implica\]\[\overthen{\recuadro{$d_\mu,d_\nu\in\mathcal L^2$}}|\dcov(\theta)|\leq\sqrt{\int d_\mu^2\dtheta^2\:\int d_\nu^2\dtheta^2}<+\infty\text{.}\]
A third approach is to derive a particular case of the AM-GM inequality (and also of Young's):
\[(d_\mu\pm d_\nu)^2\geq0\eqv\frac{d_\mu^2+d_\nu^2}{2}\geq\mp d_\mu d_\nu\eqv\frac{d_\mu^2+d_\nu^2}{2}\geq|d_\mu d_\nu|\text{,}\]
Anyhow, the key step is to show that the integrals on the right-hand side are finite. For instance, in the case of $d_\mu$: 
\[\int d_\mu(x,x')^2\dtheta^2((x,y)(x',y'))\overeq{Fubini}\iint d_\mu(x,x')\dtheta(x,y)\dtheta(x',y')\overeq{ACOV}\]
\[=\int d_\mu(x,x')\dmu^2(x,x')\overless{$d_\mu\in\mathcal L^2(\mu\times\mu)$}+\infty\text{.}\]
where the acronym ``ACOV'' stands for \emph{abstract change of variables}, which in this case takes a projection as the change of variables function. More formally, let $f$ be a measurable function in the following diagram:
\[\left(\spx\times\spy,\Borel{\spx}\tensor\Borel{\spy},\theta\right)\overarrow{$\pi_1$}\left(\spx,\Borel{\spx}\right)\overarrow{$f$}\left(\mathbb{R},\Borel{\mathbb{R}}\right)\text{.}\]
When $f\in\mathcal L^1(\theta\comp\pi_1^{-1})$, the aforementioned ACOV theorem ensures that:
\[\int_{\pi_1(\spx\times\spy)}f\wrt(\theta\comp\pi_1^{-1})=\int_{\spx\times\spy}(f\comp\pi_1)\dtheta\]
or, recalling that $\mu\overeq{def.}\theta\comp\pi_1^{-1}$:
\[\pushQED{\qed}
\int_{\spx}f(x)\dmu(x)=\int_{\spx\times\spy}f(x)\dtheta(x,y)\text{.}\qedhere
\popQED\]

The different integrability checks that have been conducted so far allow to write \emph{dcov} in terms of expected values. Taking $X\sim\mu\in M_1^1(\spx)$ and $Y\sim\nu\in M_1^1(\spy)$, with joint distribution $\theta:=\Prob\comp\binom{X}{Y}^{-1}$, their distance covariance is given by:
\[\dcov(X,Y)\overdef{Abuse}\dcov(\theta)=\E[d_\mu(X,X')d_\nu(Y,Y')]=\]
\[=\E\Big\lbrace\Big(\dx(X,X')-\E[\dx(X,X')|X]-\E[\dx(X,X')|X']+\E[\dx(X,X')]\Big)\cdot\]
\[\cdot\Big(\dy(Y,Y')-\E[\dy(Y,Y')|Y]-\E[\dy(Y,Y')|Y']+\E[\dy(Y,Y')]\Big)\Big\rbrace\text{;}\]
where primed letters refer to independent and identically distributed copies of the corresponding random element.

Finally, note that \emph{dcov} is always an association measure, in the sense that it vanishes under independence:
\[\dcov(\mu\times\nu)=\int d_\mu d_\nu\wrt(\mu\times\nu)^2\overeq{Fubini}\]
\[=\left(\int\dx\dmu^2-2\int\amu\dmu^2+\int D(\mu)\dmu^2\right)\left(\int\dy\dnu^2-2\int\anu\dnu^2+\int D(\nu)\dnu^2\right)=\]
\[=[D(\mu)-2D(\mu)+D(\mu)][D(\nu)-2D(\nu)+D(\nu)]=0\text{.}\]
Moreover, under certain conditions, \emph{dcov} is nonnegative and it can be rescaled into the interval $[0,1]$ (see~\ref{Def_dcor}), becoming a normalised association measure \citep[pages 375--376]{Bishop}.

\section{Distance covariance in negative type spaces}\label{TN}
The fact that:
\[\theta=\mu\times\nu\implica\dcov(\theta)=0\text{,}\]
makes it natural to wonder which spaces ensure that the reciprocal implication also holds. The answer is: \emph{strong negative type} spaces, since in them $\dcov(\theta)$ can be presented as an injective function of $\theta-\mu\times\nu$.

In order to explain this, negative type spaces will be firstly introduced (\S~\ref{Def_TN}), as they are the ones in which \emph{dcov} admits the aforementioned representation (although injectivity is not guaranteed). Then the strong version of this condition will be defined (\S~\ref{TNF}) and a pivotal result will be put forward: strong negative type is not only a necessary condition for \emph{dcov} to characterise independence, but it is also sufficient (with a little exception, by no means restrictive).

\subsection{Metric spaces of negative type}\label{Def_TN}
The concept of negative type is not a recent invention \citep{Wilson} and it has recently been enjoying its ``second youth'': firstly, because of its role in computational algorithmics (\citeauthor{Deza:Laurent}, \citeyear{Deza:Laurent}, \S~6.1.; \citeauthor{Naor}, \citeyear{Naor}) and, more recently, in relation to the \emph{energy of data} \citep{TEOD}.

The metric space $\xd$ is said to be of negative type if and only if:
\[\forall n\in\Z^+;\:\forall x,y\in\spx^n:\:2\sum_{i,j=1}^n\dx(x_i,y_j)\geq\sum_{i,j=1}^n[\dx(x_i,x_j)+\dx(y_i,y_j)]\text{.}\]
The analytic expression above has the following geometrical interpretation: given $n$ red points and as many blue ones, the sum of the distances among the $2n^2$ ordered pairs of the same colour is not less than the corresponding sum for different colours. Moreover, this condition can be stated in another way, that is apparently more general, which is the
\emph{conditionally negative definiteness} of the metric. However, both are actually equivalent (which can be checked by taking repetitions of the points and recalling that $\Q$ is dense in $\mathbb{R}$):
\[\forall n\in\N;\:\forall x\in\spx^n;\:\forall\alpha\in\mathbb{R}^n,\sum_{i=1}^n\alpha_i=0:\;\sum_{i,j=1}^n\alpha_i\alpha_j\dx(x_i,x_j)\leq0\text{.}\]
\setcitestyle{square}%
This is not to say that negative type metric spaces are the ones in which the metric acts like a negative definite kernel (such as the ones thoroughly studied by \citet{Klebanov} and \citet{Berg}).\setcitestyle{round}%
However, an equivalent definition in terms of the definiteness of a certain kernel exists. Namely, $\xd$ is a negative type space if and only if there is a point $o\in\spx$ so that the \emph{absolute antipodal divergence}
\[d_o(x,y):=\dx(x,o)+\dx(y,o)-\dx(x,y),\;(x,y)\in\spx^2\]
is definite positive.

There are many familiar examples of negative type spaces, like the Euclidean ones and, more generally, all Hilbert spaces (as it will be explained in~\ref{Hilbert}).

\subsection{Representation in Hilbert spaces}\label{Hilbert}
Now some results involving Hilbert spaces are to be presented. For the sake of simplicity, assume that the scalar field is $\mathbb{R}$ in every case, but, as a general rule, every statement that will be made is also true for $\mathbb{C}$, \textit{mutatis mutandi}. This can be proven by realifying or complexifying \citep[pages 132--135 of][]{Jakobsen}, according to the case.

It will be necessary to integrate functions $f:\spx\longrightarrow\h$ which have a Hilbert space as their codomain. Had $\spx$ not been assumed to be separable (see~\S~\ref{Sep}), as in \citet{Lyons}, the spaces $\h$ that arise later on would not necessarily be separable, which would only allow to perform weak integration \citep{Pettis}, and not the strong one \citep{Bochner}. Given $\mu\in M(\spx)$, if $f$ is a Pettis-integrable (or, specifically, scalarly $\mu-$integrable), the integral $I\in\h$ is unambiguously defined by its commutativity with respect to every map of the dual space $\h^*$:
\[I=\int_{\spx}f\dmu\eqv\forall h^*:\h\longrightarrow\mathbb{R}\text{ linear and continuous},\;h^*(I)=\int_{\spx}(h^*\comp f)\dmu\text{.}\]
Hereinafter, every Hilbert space that will arise is going to be separable, which means that Pettis integrals are Bochner integrals.

After these technical remarks, the \emph{Schoenberg's theorem} (\citeauthor{Schoenberg:1937}, \citeyear{Schoenberg:1937} and \citeyear{Schoenberg:1938}), can be stated. It characterises negative type spaces $\xd$ as those such that $\xsd$ can be isometrically embedded into a Hilbert space:
\[\exists\;\h\text{ Hilbert space};\:\exists\:\varphi:\spx\longrightarrow\h;\:\forall x,y\in\spx:\;\norma{\varphi(x)-\varphi(y)}_{\h}^2=\dx(x,y)\text{.}\]

For a simple proof, using the absolute antipodal divergence (see~\ref{Def_TN}), refer to \citet[theorem 3.7]{Jakobsen}, that corrects \citet{Lyons}. Regardless of this, Schoenberg's theorem ensures that the separability of the original metric spaces (\S~\ref{Sep}) is inherited by all the Hilbert spaces that arise. Before the Hilbert space representation of \emph{dcov} can be tackled, the \emph{barycentre operator} has to be defined: given an isometric map $\varphi:\xsd\longrightarrow\h_1$ (like the one on the preceding theorem) and $\mu\in M^1(\spx)$, the following Pettis integral always exists
\[\beta_\varphi(\mu):=\int_\spx\varphi\dmu\in\h_1\]
and it is called \textit{barycentre}, because it is the average of a $\h_1$-field over $\spx$ according to the distribution given by $\mu$ (thus resembling the geometrical idea of a gravity centre). In fact, if $X\sim\mu\in M_1^1(\spx)$,
\[\beta_\varphi(\mu)=\E[\varphi(X)]\text{.}\]
On the other hand, if $\psi:\ysd\longrightarrow\h_2$ is also isometric, the barycentre of the tensor product $\varphi\otimes\psi$ for $\theta\in M^{1,1}(\spx\times\spy)$ is defined as:
\[\beta_{\varphi\tensor\psi}(\theta):=\int_{\spx\times\spy}(\varphi\tensor\psi)\dtheta\in\h_1\tensor\h_2\text{.}\]

More importantly, if $(\mu,\nu)$ are the marginals of $\theta\in M_1^{1,1}(\spx\times\spy)$, the following equality holds:
\[\dcov(\theta)=4\norma{\beta_{\phi\tensor\psi}(\theta-\mu\times\nu)}_{\h_1\tensor\h_2}^2\text{.}\]
In conclusion, \emph{dcov} will characterise independence in those spaces in which the previous kernel is injective, that are going to be dealt with right below.

\subsection{Strong negative type space}\label{TNF}\setcitestyle{square}%
If $\xd$ has negative type, one can derive the following inequality (whose proof is surprisingly long \citep[lemma 3.16]{Jakobsen}):
\[\forall\mu_1,\mu_2\in M^1_1(\spx):\:D(\mu_1-\mu_2)\leq0\text{.}\]
On top of that, if the operator $D$ separates probability measures (with finite first moments) in $\xd$, that space is said to have \emph{strong} negative type:
\[D(\mu_1-\mu_2)=0\eqv\mu_1=\mu_2\text{.}\]

The extended Schoenberg's theorem shows the equivalence of the strong negative type of $\xd$ and the existence of an isometric map $\varphi:(\spx,\sqrt{\dx})\longrightarrow\h_1$ such that $\beta_\varphi$ is injective. Furthermore, for strong negative type $\spx$ and $\spy$, two isometric maps $\varphi:(\spx,\sqrt{\dx})\longrightarrow\h_1$ and $\psi:(\spy,\sqrt{\dy})\longrightarrow\h_2$ can be found so that $\beta_{\varphi\tensor\psi}:M^{1,1}(\spx\times\spy)\longrightarrow\h_1\tensor\h_2$ is injective. As a result, whenever $\spx$ and $\spy$ have strong negative type,
\[\dcov(X,Y)=0\eqv X,Y\text{ independent}\]
holds for any random element $Z=(X,Y):\Omega\longrightarrow\spx\times\spy$.
\setcitestyle{round}%

Thus, the strong negative type of marginal spaces is a \emph{sufficient} condition for the equivalence above to hold, but is it also \emph{necessary}? The answer is \emph{yes}, but with the exception of a ``pathological'' case.

If $\yd$ was not of strong negative type (symmetrically for $\spx$), it is indeed possible to find $\theta\in M^{1,1}_1(\spx\times\spy)$ so that:
\[\dcov(\theta)=0\text{ and, at the same time, }\theta\neq(\theta\comp\pi_1^{-1})\times(\theta\comp\pi_2^{-1})\text{;}\]
whenever $\min\llaves{\card\spx,\card\spy}>1$. Such $\theta$ can be constructed as follows:
\[\theta:=\frac{\delta_{x_1}\times\nu_1+\delta_{x_2}\times\nu_2}{2}\text{;}\]
where $\nu_1,\nu_2$ are two different measures in $M_1^1(\spy)$ so that $D(\nu_1-\nu_2)=0$, while $x_1,x_2\in\spx$ are two distinct points. For each $x\in\spx$, $\delta_x\in M^1(\spx)$ denotes point mass at $x$: $\delta_x(x)=1$.

This way, the aforementioned pathological case consists of one of the marginal spaces being a singleton. Such exception is not a restriction because, whenever $\card\spy=1$ (symmetrically for $\spx$), $\dcov\equiv0$ (since $d_\nu\equiv0$) and every $\theta\in M_1^{1,1}(\spx\times\spy)$ is the product of its marginals. To see this last part, note that:
\[\spy=\{y\}\implica\Borel{\spy}=\llaves{\emptyset,\{y\}}=\llaves{\emptyset,\spy}\text{.}\]
And consequently, for $B\in\Borel{\spy}$,
\[\forall A\in\Borel{\spx},\;\theta(A\times B)=\begin{cases}\theta(A\times\emptyset)=\theta(\emptyset)=0=\mu(A)\nu(\emptyset)\\
\theta(A\times\spy)=\theta\left[\pi_1^{-1}(A)\right]\equiv\mu(A)=\mu(A)\nu(\spy)\end{cases}\text{;}\]
and so $\theta=\mu\times\nu$. 
This analytical result is the formalisation of the intuitive notion that, if a random element $Y$ has constantly a certain value, the observations of any other random $X$ are bound to be independent of the ones of $Y$.

After the previous theoretical discussion, the interest of identifying practical examples of strong negative type spaces is clear. With regard to this, for the scope of the present article (and for most real data applications), it suffices to know that all separable Hilbert spaces have strong negative type. Although this is an unsurprising result, its proof is by no means straightforward \citep[pages 49--60]{Jakobsen}.

\section{Distance correlation in metric spaces}\label{dcor}

\subsection{The association measure \emph{dcor}}\label{Def_dcor}

Like previously, let $(X,Y)\sim\theta\in M_1^{1,1}(\spx\times\spy)$ have marginals $(\mu,\nu)$, where $\xd$ and $\yd$ are two separable metric spaces. Then, the following inequalities hold:
\[|\dcov(X,Y)|\leq\sqrt{\dvar(X)\dvar(Y)}\leq D(\mu)D(\nu)\text{;}\]
where $\dvar(X):=\dcov(X,X)$. If, in addition, $\xd$ and $\yd$ have negative type:
\[\dcov(X,Y)=4\norma{\beta_{\varphi\times\psi}(\theta-\mu\times\nu)}^2_{\h_1\tensor\h_2}\geq0\text{.}\]
In this context, \emph{distance correlation} (for metric spaces) is defined as:
\[\dcor(X,Y):=\frac{\dcov(X,Y)}{\sqrt{\dvar(X)\dvar(Y)}}\in[0,1]\]
whenever the denominator is nonzero. For nondegenerate cases, this will not be a matter of concern, for $\dvar(X)$ only reaches the extreme values of its range $[0,D(\mu)^2]$ when it is concentrated on one or two points (respectively):
\[\dvar(X)=0\eqv\exists\:x\in\spx,\;\mu=\delta_x\text{ ``$\mu-$almost surely'';}\]
\[\dvar(X)=D(\mu)^2\eqv\exists\:x,x'\in\spx,\;\mu=\frac{\delta_x+\delta_{x'}}{2}\text{ ``$\mu-$almost surely''.}\]
When $\dvar(X)=0$, as in the Euclidean case, $\dcor(X,Y):=0$.

\subsection{\emph{dcor} in Euclidean spaces}\label{dcor=dCor^2}
In has already been shown that \emph{dcor} has range $[0,1]$ and is zero if and only if there is independence, which recapitulates the property for Euclidean spaces (\S~\ref{SRB}). Indeed, it is possible to prove (via the Hilbert space representations introduced in~\ref{Hilbert}) that, when $\xd$ and $\yd$ are (finitely dimensional) Euclidean spaces, the notion of distance correlation of \S~\ref{Def_dcor} \citep{Lyons} generalises the square of the one in \S~\ref{SRB} \citep{SRB}:
\[\dcov(X,Y)=\dCov(X,Y)^2;\;\dcor(X,Y)=\dCor(X,Y)^2\text{.}\]
For $\theta\in M_1^{2,2}(\spx\times\spy)$, $\dcov(X,Y)$ becomes a product of expectations. By expanding it and simplifying, one can easily get the generalisation of Brownian distance covariance \citep[theorems 7--8]{Brownian} to general metric spaces:
\[\dcov(X,Y)\!=\!\E[\dx(X,X')\dy(Y,Y')]+\E[\dx(X,X')]\E[\dy(Y,Y')]-\]\[-2\E[\dx(X,X')\dy(Y,Y'')]\text{.}\]
In conclusion, \emph{dcov} satisfactorily extends \emph{dCov} squared.

\section{Nonparametric test of independence in metric spaces}\label{Test dcor}

\subsection{Kernel associated to \emph{dcov}}\label{h}
The following map will be key to the construction of the sample version of \emph{dcov}:
\begin{align*}
h:\;&(\spx\times\spy)^6\flecha\mathbb{R}\\
&\big((x_i,y_i)\big)_{i=1}^6\flechita f_{\spx}(x_1,x_2,x_3,x_4)\,f_{\spy}(y_1,y_2,y_5,y_6)\text{;}
\end{align*}
where, for $\spz\in\llaves{\spx,\spy}$,
\[f_{\spz}(z):=\dz(z_1,z_2)+\dz(z_3,z_4)-\dz(z_1,z_3)-\dz(z_2,z_4),\:z\in\spz^4\text{.}\]
\setcitestyle{square}%
The functions $f_{\spz}$ and $h$ are clearly measurable and proving their integrability can be accomplished by sequentially deriving inequalities from the triangle inequality (see pages 148--150 of \citet{Jakobsen} for the correction of the attempt by \citet{Lyons}). Integrating these functions is pretty straightforward. Firstly, for $f_{\spx}$:\setcitestyle{round}%
\[\int_{(\xy)^2}f_{\spx}(x_1,x_2,x_3,x_4)\dtheta^2((x_3,y_3),(x_4,y_4))\overeq{ACOV}\]
\[=\dx(x_1,x_2)-\amu(x_1)-\amu(x_2)+D(\mu)\equiv d_\mu(x_1,x_2),\;(x_1,x_2)\in\spx^2\text{;}\]
where $\theta\in M_1^{1,1}(\xy)$ has marginals $(\mu,\nu)$. Given that the same (\textit{mutatis mutandi}) holds for $f_{\spy}$,
\[\dcov(\theta)=\int_{(\xy)^2}d_\mu(x_1,x_2)d_\nu(y_1,y_2)\dtheta^2\left((x_1,y_1),(x_2,y_2)\right)=\int_{(\xy)^6}h\dtheta^6\text{.}\]
This means that, if $(X_i,Y_i)_{i=1}^6$ denotes a vector that contains $6$ random elements that are independent and identically distributed to $(X,Y)\sim\theta$,
\[\dcov(\theta)=\E\left[h\left((X_i,Y_i)_{i=1}^6\right)\right]\]
and, consequently, its sample version is a $V-$statistic, as the ones that \citet{Lyons} derived (erroneously), as it will be shown next.

\subsection{Empirical distance covariance}\label{UV}

For $n\in\mathbb{Z}^{+}$, the following notation will be used for the \emph{empirical measure} associated to a certain sample $\llaves{(X_i,Y_i)}_{i=1}^n\iid(X,Y)\sim\theta$: 
\[\theta_n:=\frac{1}{n}\sum_{i=1}^n\delta_{(X_i,Y_i)}:\Omega\longrightarrow M_1^{1,1}(\xy)\text{.}\]
A few routine computations yield that the natural estimator
\[\dcovh(\theta):=\dcov(\theta_n)\]
is, unsurprisingly, the $V-$statistic with (nonsymmetric) kernel $h$:
\[\dcov(\theta_n)=\frac{1}{n^6}\sum_{i_1=1}^n\cdots\sum_{i_6=1}^n h\left((X_{i_\lambda},Y_{i_\lambda})_{\lambda=1}^6\right)\equiv V_n^6(h)\text{.}\]
On the other hand, it is logical to consider the analogous $U-$statistic as an alternative estimator, which will be shown to require less stringent conditions to behave satisfactorily than $\dcov(\theta_n)$. For $n\geq7$, let:
\[\tilde U_n^6(h):=\frac{1}{6!\binom{n}{6}}\sum_{\llaves{i_\lambda}_{\lambda}\subset[1,n]\cap\mathbb{Z}\text{ different}}h\left((X_{i_\lambda},Y_{i_\lambda})_{\lambda=1}^6\right)\text{;}\]
where the tilde indicates that this is not a $U-$statistic \textit{sensu stricto}, but rather one built upon a kernel that is nonsymmetric. To correct this, let $\bar h$ be the symmetrisation of $h$:
\[\bar h(z):=\frac{1}{6!}\sum_{\sigma\in S_6}h\left(z_{\sigma(j)}\right)_{j=1}^6\equiv\frac{1}{6!}\sum_{\sigma\in S_6}h(z_\sigma),\:z\in(\xy)^6\text{;}\]
where $S_6:=\llaves{\sigma:[1,6]\cap\mathbb{Z}\flecha[1,6]\cap\mathbb{Z}:\;\sigma\text{ bijective}}$ is the symmetric group of order $6$.
So $\tilde U_n^6(h)$ is the $U-$statistic based on $\bar h$:
\[\tilde U_n^6(h)=\frac{1}{\binom{n}{6}}\sum_{i_1<\ldots<i_6}\bar h\left((X_{i_\lambda},Y_{i_\lambda})_{\lambda=1}^6\right)\text{.}\]
The analogous for the $V-$statistic also holds:
\[\forall\:\sigma\in S_6,\;\dcov(\theta_n)\equiv V_n^6(h)=\int_{(\xy)^6}h(z)\dtheta_n^6(z)\overeq{Fubini}\]
\[=\int_{(\xy)^6}h(z)\dtheta_n^6(z_{\sigma^{-1}})\overeq{ACOV}\int_{(\xy)^6}h(z_\sigma)\dtheta_n^6(z)=V_n^6(\bar h)\]
and the same arguments can prove that $\dcov(\theta)=\int\bar h\dtheta^6$.

Now that the usual symmetric kernels can be used, it is possible to resort to the \emph{strong law of large numbers} (SLLN) for $U-$statistics \citep{Hoeffding} to infer that, for $\theta\in M_1^{1,1}(\xy)$,
\[\tilde U_n^6(h)\aslim\dcov(\theta)\text{.}\]
\begin{sloppypar}
\citet{Lyons} mistook the hypotheses of the aforementioned Hoeffding theorem for the ones of the SLLN for $V-$statistics \citep[page 274]{Gine:Zinn}. The weakest conditions under which the SLLN for $V-$statistics hold in this context are: ${\theta\in M_1^{5/3,5/3}(\xy)}$ \citep[theorem 5.5]{Jakobsen}. In other words, the finiteness of moments of order $\frac{5}{3}$ suffices to ensure asymptotic consistency:
\end{sloppypar}
\[V_n^6(h)\aslim\dcov(\theta)\text{.}\]

\subsection{Null distribution of the test statistic}\label{NullDist}
If $\theta\in M_1^{1,1}(\xy)$ is the product of its marginals $(\mu,\nu)$ and these are nondegenerate, the asymptotic distributions of the estimators introduced in~\ref{UV} are:
\begin{align*}
nV_n^6(h)&\distrilim\sum_{i=1}^\infty\lambda_i(Z_i^2-1)+D(\mu)D(\nu)\text{;}\\
n\tilde U_n^6(h)&\distrilim\sum_{i=1}^\infty\lambda_i(Z_i^2-1)\text{;}
\end{align*}
where $\llaves{Z_i}_{i\in\Nstar}\iid\Normal(0,1)$ and where $\llaves{\lambda_i}_{i\in\Nstar}$ are the eigenvalues (with multiplicity) of the linear operator $S:\mathcal L^2(\theta)\flecha\mathcal L^2(\theta)$ that maps $f$ into $S(f):\xy\flecha\mathbb{R}$, which is defined as:
\[S(f)(x,y):=\int_{\xy}d_\mu(x,x')d_\nu(y,y')f(x',y')\dtheta(x',y'),\:(x,y)\in\xy\text{.}\]

The original attempt of proving the result for the $V-$statistic \citep{Lyons} included some incorrect arguments to conclude that $\sum_{i=1}^\infty\lambda_i=D(\mu)D(\nu)$. \citet{Ups} states that the previous identity does hold as long as both marginal spaces have negative type, but the justification of this is somewhat abstruse. In case of it being true, it would be the exact same asymptotic distribution that \citet{SRB} had derived. Anyhow, this cannot be brought to practical usefulness (as in~\ref{SRB}), since the eigenvalues $\llaves{\lambda_i}_i$ depend on $\theta$ (unknown) and cannot be easily estimated. The most logical approach to this is, once again as in~\ref{SRB}, a resampling strategy. One way of arguing for this procedure would be to summon the results of \citet{Arcones:Gine}, that ensure that approximating the thresholds for the test statistic via na\"ive bootstrap leads to a consistent resampling technique, as $\bar h$ satisfies the integrability condition required by those authors.

\end{document}